Dmytro Taranovsky
February 26, 2012
Small change: December 12, 2016

# Determinacy and Fast-growing Sequences of Turing Degrees

**Abstract:** We discuss sufficiently fast-growing sequences of Turing degrees. The key result is that, assuming sufficient determinacy, if $\varphi$ is a formula with one free variable, and $S$ and $T$ are sufficiently fast-growing sequences of Turing degrees of length $\omega_1$, then $\varphi(S) \Leftrightarrow \varphi(T)$. In the second part (below), we define degrees for subsets of $\omega_1$ analogous to Turing degrees, and prove that under sufficient determinacy and CH, all sufficiently high degrees are also effectively indistinguishable.

## Fast-growing Sequences of Turing Degrees

A consequence of definable determinacy is that for every definable property, either it or its negation holds on a cone of Turing degrees. In other words, if $T$ is a sufficiently high Turing degree, then its definable properties are independent of $T$. Moreover, if $T_1$ is a sufficiently high Turing degree and $T_2$ is sufficiently high above $T_1$, then definable properties of $(T_1, T_2)$ are independent of $(T_1, T_2)$. Here, "definable" depends on determinacy assumed; for projective determinacy "definable" is "definable in second order arithmetic". (Recall that projective determinacy asserts that every two-player perfect information game of length $\omega$ on integers (that is with moves being integers) and projective payoff set is determined, that is one of the players has a winning strategy.) We can ask whether the analogous property holds for infinite sufficiently fast-growing sequences of Turing degrees, and the answer turns out to be yes in a very strong way.

**Theorem 1** (ZFC + determinacy for games on integers of length $\omega_1$ and ordinal definable payoff):
Let $\delta > \omega_1$ be an ordinal and $r$ a real number. There is a function $f$ from countable sequences of Turing degrees to Turing degrees such that for all formulas $\varphi$ with two free variables, for every sequence of Turing degrees $S$ of length $\omega_1$ with
$\forall \alpha < \omega_1 \, f(S|\alpha) \leq S(\alpha)$, the truth of $\varphi(r, S)$ in $V_\delta$ is independent of $S$.

**Proof:** Determinacy in the theorem also implies determinacy for payoff ordinal definable from a real. Pick a formula $\varphi$ and consider the game of length $\omega_1$ with the payoff set $\{x : \varphi_\delta(r, S(x))\}$ where $\varphi_\delta$ is relativization of $\varphi$ to $V_\delta$ and $S(x)(\alpha)$ is the Turing degree of $[x_{\omega \cdot \alpha}, x_{\omega \cdot \alpha + \omega})$ ($x$ is a sequence of integers of length $\omega_1$). By determinacy, one of the players has a winning strategy. Because the payoff depends only on the Turing degrees (and with the help of a well-ordering of real numbers), there is a winning strategy that depends only on the last finite set of moves and the Turing degrees of the previous ω-sequences of moves. Pick such a strategy, and let $f_\varphi$ be such that $f_\varphi(S|\alpha)$ is the Turing degree for the strategy for moves $\omega \cdot \alpha$ up to $\omega \cdot \alpha + \omega$ where $S|\alpha$ is the sequence of Turing degrees generated by the moves before $\omega \cdot \alpha$. Since the other player can force the resulting sequence of Turing degrees to be any $S$ with $f_\varphi(S|\alpha) \leq S(\alpha)$, the truth of $\varphi(r,S)$ in $V_\delta$ is independent of $S$ for all such $S$. Finally, pick $f$ with $\forall S \in w^{<\omega_1} \forall \hat\varphi\, f(S) \geq f_\varphi(S)$, and $f$ satisfies the theorem.

**Notes:**

- By varying $\delta$, we can code-in an arbitrary ordinal as a parameter. Also, $\delta > \omega_1$ so that $\omega_1$ (and hence $S$) is in $V_\delta$.
- The determinacy assumption in the theorem is consistent relative to a Woodin limit of Woodin cardinals and (but that should be unnecessary) a measurable cardinal above them (see for example *The determinacy of long games* reference). A leisurely introduction to a similar determinacy hypothesis can be found in the "Determinacy Maximum" reference.
- "Ordinal definable payoff" can be replaced with another sufficiently robust notion of definability such as definability in third order arithmetic provided that $\varphi$ is restricted to that class of definability. Also, if we allow any payoff ordinal definable from a countable sequence of ordinals, then we can allow $r$ in $\varphi$ to be any countable sequence of ordinals.
- Weaker forms of determinacy translate into weaker forms of the theorem. For example, the theorem holds with $\omega_1$ replaced by a countable ordinal α and determinacy for games of length $\omega \cdot \alpha$ and OD(R) payoff.
- Projective determinacy suffices for the theorem restricted to the set of all $\varphi$ expressible in second order arithmetic and finite sequences of Turing degrees (of the same length). For that set of $\varphi$, the minimum notion of sufficiently high that suffices for finite sequences is: $r_2$ has a sufficiently high Turing degree above that of $r_1$ iff for every natural number $n$, $\Sigma^1_n(r_1)$ truth is computable from $r_2$.
- The theorem holds even if Turing degrees are replaced by elementary time degrees: $X$ is elementary time computable from $Y$ iff there is $n$ such that $X$ can be computed from $Y$ using time bounded by a stack of $n$ exponentials. However, using polynomial time degrees would not work since coding the opponent's strategy into the play is exponential (for games on {0,1} of length $\omega$); relativized P=NP toggles even for arbitrarily high polynomial time degrees.

**Definition:** We say that $\varphi$ holds for every sufficiently fast-growing sequence of Turing degrees (of a particular length $\leq \omega_1$) if there is $f$ (from countable sequences of Turing degrees to Turing degrees) such that $\varphi$ holds for all $S$ with $\forall \alpha < |S|\, f(S|\alpha) \leq S(\alpha)$ ($|S|$ is the order type of $S$).

For $|S| > \omega$, the notion of sufficiently fast-growing is not closed under subsequences. Every infinite sequence has an uncountable number of subsequences, so for example, for every $S$, there are subsequences of $S|\omega$ that are not recursive in $S(\omega)$. However, the notion of 'sufficiently fast-growing' is well-behaved in that for every countable set of conditions (each represented by $f_n$), there is a sequence that satisfies all of them.

In general, a notion of sufficiently fast-growing sequences of Turing degrees is a function that maps each countable sequence of Turing degrees into a nonempty upwards-closed set of Turing degrees (which are the permissible degrees for the next element of the sequence).

**Proposition 2:** Given a notion of sufficiently fast-growing sequences of Turing degrees, there is a stricter notion $f$ such that (1) for successor $|S|$, $f(S)$ depends only on on the last element of $S$ provided that $\forall \alpha < |S|\, S(\alpha) \in f(S|\alpha)$, and (2) (strong monotonicity)

if $T \in L_\delta(S)$ for the least $\delta$ for which $L_\delta(S)$ satisfies ZF\P, then $f(S) \subset f(T)$.
**Proof:** Let $f_1$ be the initial notion strengthened so that each element of $f_1(S)$ can enumerate $S$, which reduces (1) to (2). $f$ defined by $f(S) = \cap f_1(T)$, where $T$ is obtained from $S$ by zero or more (actually, exactly one suffices) applications of (2), satisfies the proposition.

For reference, $\text{HOD}_B^A(C)$ consists of hereditarily ordinal definable sets as computed in $A$ (or if $A$ is omitted, in $V$) allowing $B$, $C$, and elements of $C$ as parameters.

**Corollary 3 (to Theorem 1)** (ZFC + determinacy for games on integers of length $\omega_1$ and ordinal definable payoff):
   *(a)* If $S$ is a sufficiently fast-growing sequence of Turing degrees of length $\omega_1$, then all reals in $\text{HOD}_S$ are ordinal definable, and hence $\mathbb{R} \cap \text{HOD}_S$ is countable.
   *(b)* If $\delta$ is a countable ordinal and (relative to $\delta$) $S$ is a sufficiently fast-growing sequence of Turing degrees of length $\omega_1$, then $\text{HOD} \cap V_\delta = \text{HOD}_S \cap V_\delta$.

**Proof:** *(a)* Assume contrary, and let $\alpha$ be the least ordinal such that (under the canonical well-ordering of $\text{HOD}_S$), $\alpha$-th real is not in $\text{HOD}$, and let $r(S)$ denote that real. From the theorem, $r(S)$ is independent from $S$, which contradicts $r(S)$ not being ordinal definable.
   *(b)* Assume contrary, and let $W$ be a canonical well-ordering of $\text{HOD}_S$ (mapping ordinals to sets) such that sets of lower rank come before sets of higher rank and $W$ agrees with a canonical well-ordering of $\text{HOD}$ up to the least set not in $\text{HOD}$. Let $\alpha$ be the least ordinal such that $W(\alpha)$ is in $\text{HOD}_S \cap V_\delta$ but not in $\text{HOD}$. Under the determinacy assumption, $\omega_1$ is inaccessible (and even measurable) in $\text{HOD}$. Because $\alpha \leq f(\delta)$ for some definable function $f$ (specifically, the order of $\text{HOD} \cap V_\delta$ under the canonical well-ordering of $\text{HOD}$) from countable ordinals to countable ordinals, $\alpha$ is independent of $S$ (to prove, apply the theorem to each $\alpha \leq f(\delta)$ and take upper bound for the rate of growth). $W(\alpha)$ must depend on $S$, so there is the least ordinal $\beta$ (with $\beta \leq f(\delta)$) such that the truth of $W(\beta) \in W(\alpha)$ depends on $S$, which contradicts the theorem.

The theorem is surprising since the freedom to pick $\omega_1$ parameters would ordinarily allow coding in an arbitrary subset of $\omega_1$. It is also surprising that even though $S$ is a sequence of $\omega_1$ disjoint countable sets of reals, the set of reals ordinal definable from $S$ is countable.

It seems reasonable that the least disagreement between $\text{HOD}$ and $\text{HOD}_S$ corresponds to a measurable cardinal in $\text{HOD}_S$ and an associated elementary embedding. For each s∈S, the least definable from S ordinal corresponding to *s* appears to be the supremum of *s*-recursive well-orderings, and a reasonable conjecture is that it is measurable in $\text{HOD}_S$. Also, based on results in "Large Cardinals from Determinacy" chapter of the Handbook of Set Theory, it seems reasonable that each element of $S$ leads to a Woodin cardinal in $\text{HOD}_S$.

An interesting project would be to work out the theory of definability from a sufficiently fast-growing sequence of Turing degrees $S$. Here are some partial results.

**Proposition 4:**
   *(a)* If $S$ is a sufficiently fast-growing sequence of $n < \omega$ Turing degrees, then

$\text{HOD}_S^{L(S)}$ contains all $\Delta^1_{n+2}$ (and $\Sigma^1_{n+2}$) reals, and assuming projective determinacy (or just absence of a projective sequence of $\omega_1$ distinct reals), all reals that are $\Delta^1_{n+2}$ in a countable ordinal.

(b) Let $\delta$ be an ordinal. If $S$ is a sufficiently fast-growing sequence of $\omega$ Turing degrees, then $\text{HOD}_S^{L(S)}$ contains (the real number coding) the truth predicate of $L_\delta(\mathbb{R})$. If $\text{HOD}^{L(\mathbb{R})} \cap \mathbb{R}$ is countable (for example if the axiom of determinacy holds in $L(\mathbb{R})$), then $S$ can be chosen independent of $\delta$.

**Proof:** *(a)* $S$ can be coded by a real number, and provably in ZFC and uniformly in the code for $S$, $S$ can be used to convert $n$ real quantifiers into $n$ arithmetic quantifiers. For example, $\exists X \in \mathbb{R}\, \forall Y \in \mathbb{R}\, \varphi(X, Y) \Leftrightarrow \exists X \leq_T r_0\, \forall Y \leq_T r_1\, \varphi(X, Y)$ ($\varphi$ has no other free variables; $r_i$ can be any real number with Turing degree $S_i$). Combined with $\Sigma^1_2(S)$ correctness of $L(S)$, this allows $\text{HOD}_S^{L(S)}$ to compute $\Sigma^1_{n+2}$ truth and hence include all $\Delta^1_{n+2}$ (and $\Sigma^1_{n+2}$) reals. Furthermore, the proof relativizes to any particular countable ordinal (and appropriate $S$). Under projective determinacy, only countably many reals are $\Delta^1_{n+2}$ in a countable ordinal, so a sufficiently fast-growing $S$ will cover all of them.

*(b)* Let $S$ be such that for every $n < \omega$, there is an elementary substructure of $L_\delta(\mathbb{R})$ containing $S_n$ and with all reals having a lower Turing degree than $S_{n+1}$, and let M be the union of an $\omega$ chain of these substructures. Clearly, $\mathbb{R}^M = \mathbb{R}_S$ where $\mathbb{R}_S$ is the closure of $\bigcup S$ under Turing reducibility, so $M$ condenses to some $L_\lambda(\mathbb{R}_S)$ and hence its truth predicate is in $\text{HOD}_S^{L(S)}$. Finally, if $\text{HOD}^{L(\mathbb{R})} \cap \mathbb{R}$ is countable, then by existence of upper bounds for countable sets of notions of sufficiently fast-growing, we can pick $S$ independent of $\delta$.

**Note:** Proposition 4 and Corollary 3 can also be relativized on top of an arbitrary real $r$ (using $\text{HOD}_{r,S}$ and $S$ dependent on $r$).

**Proposition 5** (ZFC + determinacy for games on integers of length $\omega_1$ and ordinal definable payoff):

*(a)* Let $S$ be a sufficiently fast-growing sequence of Turing degrees of length $\omega_1$ and $T$ an initial segment of $S$ of limit length. Then $\mathbb{R}^{\text{HOD}_S(T)} = \mathbb{R}_T$ where $\mathbb{R}_T$ is the closure of $\cup T$ under Turing reducibility. Furthermore, there is an elementary embedding of $L(\mathbb{R}_T)$ into $L(\mathbb{R})$.

*(b)* $\text{HOD}(\mathbb{R}_T)$ agrees with $\text{HOD}(\mathbb{R})$ about statements in third-order arithmetic with parameters in $\mathbb{R}_T$ (and hence it satisfies AD).

**Proof:** *(a)* Every real in $\text{HOD}_S(T)$ is ordinal definable from $r$ and $S$ for some $r$ belonging to some T($\alpha$). Above $\alpha$, $S$ is sufficiently fast-growing relative to $r$, so all reals ordinal definable from $r$ and $S$ have Turing degree below $S_{\alpha+1}$. Because $|T|$ is a limit ordinal, it follows that $\mathbb{R}^{\text{HOD}_S(T)} = \mathbb{R}_T$.

For each $T(\alpha)$, there is a canonical elementary substructure of $L(\mathbb{R})$ containing $T(\alpha)$ and with all reals being in $T(\alpha + 1)$. The union of these substructures condenses/collapses to $L(\mathbb{R}_T)$, and the inverse of the collapse is the desired elementary embedding. Also, while $L(\mathbb{R})$ is a proper class, $\mathbb{R}^\#$ exists, so the notions of elementary substructures, etc. are well-defined.

(b) To avoid issues with undefinability of truth, work in $V_\delta$ where $\delta$ is the least ordinal such that $V_\delta \prec_{\Sigma_2} V$. Let $M$ consist of sets in $\mathrm{HOD}(\mathbb{R})$ that are definable from an element of $\mathbb{R}_T$, and $N$ consist of sets in $\mathrm{HOD}(\mathbb{R}_T)$ that are definable from $\mathbb{R}_T$ and an element of $\mathbb{R}_T$. $N$ is an elementary substructure of $\mathrm{HOD}(\mathbb{R}_T)$, and by the closure of $\mathbb{R}_T$, $M$ is an elementary substructure of $\mathrm{HOD}(\mathbb{R})$. Let $M_1$ be the transitive collapse of $M$. $\mathbb{R}^{M_1} = \mathbb{R}_T$ and $\mathrm{P}(\mathbb{R})^{M_1} \subset \mathrm{P}(\mathbb{R})^N$, so to complete the proof it suffices to show that every set of reals in $N$ is in $M_1$. Every set of reals in $N$ has the form $\{x : \varphi(r, x)\}$ for some formula $\varphi$ (with two free variables) and some $r$ in $\mathbb{R}_T$. Using Theorem 1, we can show that the truth of $\varphi(r, x)$ is independent of $T$ (for all sufficiently fast-growing $T$ of the same limit length with $r$ and $x$ in $\mathbb{R}_T$), so the set is the intersection of $\mathbb{R}_T$ with a set of reals definable from $r$, and hence is present in $M_1$.

One consequence of the proposition is that (under the determinacy assumption), for limit $|T|$ that are not sufficiently closed (including $\omega$), in $L(T)$ (or $\mathrm{HOD}_S(T)$), the continuum is a countable union of countable sets. For all limit $|T|$, $\mathrm{HOD}_S(T)$ (and hence $L(T)$) satisfies "there are no uncountable sequences of distinct reals". Also, as before, restricted forms of determinacy translate into restricted forms of the proposition.

We end with some conjectures.
**Conjectures:** For finite length $n$ (with $S$ being a sufficiently fast-growing sequence of Turing degrees of length $n$), the reals in $\mathrm{HOD}_S^{L(S)}$ are the ones $\Delta^1_{n+2}$ in a countable ordinal, equivalently, those that are in the minimal iterable inner model with $n$ Woodin cardinals. If $\alpha = |S|$ is countable in $M_\alpha$ (the minimal iterable inner model with $\alpha$ Woodin cardinals), then the reals in $\mathrm{HOD}_S^{L(S)}$ are those in $M_\alpha$, and perhaps, $\mathrm{HOD}_S^{L(S)}$ itself is an iterate of $M_\alpha$ (or an iterate of $M_\alpha$ augmented with certain iteration strategies).
We can also consider definability in $L(\mathbb{R})$. One conjecture is that for $\alpha = |S|$ countable in $M_\alpha$, the reals in $\mathrm{HOD}_S^{L(\mathbb{R})}$ are precisely those in $M_{\alpha+\omega}$. Finally, an interesting question is the nature of $L(S)$ and $\mathrm{HOD}_S^{L(S)}$ for $|S| = \omega_1$.

# Degrees of Subsets of $\omega_1$

For subsets of $\omega_1$, we can define several notions of degrees:

- Constructible reducibility: $X \leq Y$ iff $X \in L[Y]$.
- Lightface projective reducibility: $X \leq Y$ iff $X$ is definable in second order arithmetic augmented with a predicate for $Y$. Boldface projective reducibility is the same, except that it allows a real number as a parameter.
- Lightface recursive reducibility: $X \leq Y$ iff $X$ is $\Delta_1$ in $Y$ in the following sense: There are formulas $\varphi$ and $\psi$ with one free variable such that
  $\forall \alpha < \omega_1(\alpha \in X \Leftrightarrow \exists \beta < \omega_1((L_\beta[Y], \in, Y) \models \varphi(\alpha)) \Leftrightarrow \neg \exists \beta < \omega_1((L_\beta[Y], \in, Y) \models \psi(\alpha)))$. Boldface recursive reducibility is the same, except that it allows a countable ordinal as a parameter.

For each of these notions of reducibility, there is a corresponding degree notion generated by equivalence classes of subsets of $\omega_1$. Boldface recursive reducibility is finer than constructible reducibility and boldface projective reducibility. Under CH, sets of reals can be coded as subsets of $\omega_1$, but constructible reducibility (using $L(Y)$ in place of $L[Y]$) and both of projective reducibilities make sense for sets of reals directly.

Clearly, lightface recursive reducibility is the strictest notion, but it is sufficiently robust for our purposes. It turns out that under appropriate assumptions (specifically, CH and sufficient determinacy), all sufficiently high lightface recursive reducibility degrees are effectively indistinguishable.

**Theorem 6** (ZFC + determinacy for games on integers of length $\omega_1$ and ordinal definable payoff + Continuum Hypothesis):
Let $\delta > \omega_1$ be an ordinal and $r$ a real number. There is a degree $S$ for lightface recursive reducibility of subsets of $\omega_1$ such that for all formulas $\varphi$ with two free variables and all $T \geq S$, $V_\delta$ satisfies $\varphi(r, S) \Leftrightarrow \varphi(r, T)$. Furthermore, for boldface recursive reducibility, $S$ can be chosen independent of $r$.

**Proof:** Consider a game on integers of length $\omega_1$ and payoff $\{X : \varphi(r, \mathrm{degree}(X))\}$ as computed in $V_\delta$. If one player has a winning strategy, then by CH, it can be coded into a subset of $\omega_1$. Let $S$ be the degree of one such code (under any reasonable coding). Because application of the code for a strategy to a partial play is recursive in the required sense, the other player can force the degree of the play to be any degree $T \geq S$, so $\varphi(r, S) \Leftrightarrow \varphi(r, T)$. Because every countable set of of degrees has an upper bound, we can pick a degree that works for all $\varphi$. For boldface recursive reducibility, every set of $\omega_1$ degrees has an upper bound, so we can a pick a degree that works for all $\varphi$ and $r$.

**Notes:**
- Weaker forms of determinacy for games of length $\omega_1$ translate into weaker forms of the theorem (corresponding to restricted classes of φ).
- The Continuum Hypothesis is needed to be able to code the opponent's strategy into the play. We do not know what happens if CH fails.
- The theorem also applies to coarser degrees like constructible or projective reducibility.
- We do not know whether the conclusion of the theorem holds for sufficiently fast-growing pairs or sequences of degrees of subsets of $\omega_1$. We also do not know what happens (alternatively, what consistently can happen) for degrees of subsets of larger ordinals.